\documentclass[12pt,a4paper,twoside]{article}
\usepackage[top=1in, bottom=1in, right=.7in, left=1in]{geometry}
\usepackage{mathptmx}
\usepackage{amsmath}
\allowdisplaybreaks
\usepackage{float}
\newenvironment{rcases}
{\left.\begin{aligned}}
	{\end{aligned}\right\rbrace}
\usepackage{amssymb}
\usepackage{amsfonts}
\usepackage{setspace}
\usepackage{enumitem}
\usepackage[T1]{fontenc}
\usepackage{booktabs}
\usepackage{siunitx}
\usepackage{multirow}
\everymath{\displaystyle}
\usepackage[utf8]{inputenc}
\usepackage[labelfont=bf]{caption}
\usepackage{subfigure}
\usepackage{graphicx}
\usepackage[numbers]{natbib}
\usepackage{calc}
\usepackage{blindtext}
\usepackage[none]{hyphenat}
\usepackage{fancyhdr}

\usepackage[breaklinks]{hyperref}
\hypersetup{colorlinks, citecolor= green, linkcolor=blue}
\usepackage[compact]{titlesec}
\titlespacing{\section}{0pt}{*0}{*0}
\titlespacing{\subsection}{0pt}{*0}{*0}
\titlespacing{\subsubsection}{0pt}{*0}{*0}
\begin{document}

\newpage
\pagenumbering{arabic}

\begin{titlepage}

\newcommand{\HRule}{\rule{\linewidth}{0.6mm}}
	
\begin{center} 
	
	\textsc{\huge{ Impact of radiation and slip conditions on MHD flow of nanofluid past an exponentially stretched surface}}\\[1cm]
\end{center}	 
\date{}
\textsc{Diksha Sharma}$^{1}$ \textsc{Shilpa Sood}$^{2}$\\
$^{1}$ Department of Mathematics and Statistics, Career Point University Hamirpur, Himachal Pradesh.\\	
diksha.math@cpuh.edu.in\\
$^{2}$ Department of Mathematics and Statistics, Career Point University Hamirpur, Himachal Pradesh.\\	
shilpa.math@cpuh.edu.in \\\\
Corresponding Author: Diksha Sharma\\
Affiliation: Research Scholar\\
Email Address: diksha.math@cpuh.edu.in

\begin{abstract}
The current research establishes magnetohydrodynamics (MHD) boundary layer flow with heat and mass transfer of a nanofluid over an exponentially extending sheet embedded in a porous medium. During this exploration, nanoparticles, single-wall carbon nanotubes (SWCNTs) and multi-wall carbon nanotubes (MWCNTs) are recruited, while lamp fuel oil is being utilised as a base fluid for the diffusion of nano materials. The effects of warm radiation and an inclined magnetic field are included. In addition, rather than no-slip assumptions at the surface, velocity slides as well as thermal upsurge are incorporated in this study. Similarity transformations are implemented to adapt a set of partial differential equations into a system of non-linear ordinary differential equations. The bvp4c solver and Keller-box approach are employed to tackle nonlinear ordinary differential equations numerically. The significance of prominent parameters such as the Darcy-Forchheimer model, magnetic field, radiation, suction, velocity slip, and temperature jump is visually probed and addressed in depth. In fact, the evolution of the coefficient of skin friction and percentage of heat shipping (Nusselt number) for both SWCNTs and MWCNTs is presented in tabular form. The temperature goes up as the magnetic parameter rises. Temperature has been seen to be decreased as the thermal slip parameter is improved. The results indicate that SWCNTs yield a higher coefficient of skin friction and speed of heat transformation than MWCNTs.\\

\end{abstract}

\section{Introduction}

Due to the growing need for heating and cooling in industrial processes during the last few decades, researchers and engineers have focused on improving thermal conduction rate. Enhancing the thermal physical phenomena of ancient fluids with low thermal conductivity, similar to ethylene glycol, water, and oil, is a crucial task for researchers. Numerous researchers have looked into boosting the thermal resistance in a broad range of industrial implementations, along with nuclear reactors. Choi and Estman \cite{1} were the initial ones to originate the term "Nanofluid". A nanofluid is made up of nanoparticles (1-100 nm) blended together using a base liquid. Some of these nanoparticles are carbides or carbon nanotubes, oxides, metals, etc. Nanoparticles are subdivided into numerous groups based upon their size, frame, physical or chemical attributes, and a myriad of other aspects. Metallic nanoparticles and non-metallic nanoparticles are the two types of nanoparticles that exist. Metallic nanoparticles include carbides, alumina, copper, nitrides, and metal oxides, whereas non-metallic nanoparticles include carbon nanotubes and graphite. A study group at Argonne National Laboratory was the first to discuss the usage of various nanoparticles in different base fluids. They concluded that nanoparticles appear in a variety of shapes, such as spherical, cylindrical, brick-like, etc. Based on the complicated shapes of nanoparticles, the rate of heat transmission of nanofluids also alters. Carbon nanotubes, which are cylindrical-like structures whose diameters range from 1 to 50 nm, have been proven to have a better heat transmission rate than spherical or concrete block nanoparticles.\

Carbon nanotubes have a cylindrical frame and are wrapped in one or perhaps more graphene sheet layers. CNTs have a 1–2 nm diameter extended tubular structure \cite{2}-\cite{4}. The best carbon nanotubes, according to Choi and Zhang \cite{5}, are those that have been coiled up into a hexagonal arrangement of carbon molecules to construct a steam-free hollow cylinder. Sumio Lijima \cite{6} a Japanese physicist, was the first to discover CNTs in 1991 for multiple or several wall nanotubes. Carbon nanotubes (CNTs) are distributed into two different categories, i.e., single-wall carbon nanotubes (SWCNTs) and multi-wall carbon nanotubes (MWCNTs). CNTs are characterised based on the intertwined configuration of carbon molecules. Carbon nanotubes have better heat transmission properties than nanomaterials of the same volume concentration \cite{7}-\cite{8}. Carbon nanotubes (CNTs) play an imperative character in transistors, batteries, superconductors, data storage, electromagnetic shielding, semiconductors, biosensors, temperature insulators from solar storage, etc. Because of their superior chemical suitability towards biological molecules such as proteins and DNA, carbon nanotubes are additionally exploited in pharmaceuticals, healthcare products, and diagnostic devices \cite{9}. There's no way to refute the astonishing aspects of nanomaterials.\ 

Single-walled carbon nanotubes and multi-walled carbon nanotubes have many attributes, however there are several notable differences. These adaptations and frameworks features significant aesthetic and microelectronic topographies, extraordinary asset and flexibility, and superior thermal steadiness. Said et al. \cite{10} have investigated how SWCNTs nanoliquids advance the thermal recital of solar collectors instead of water. Xue et al. \cite{11} analyzed the experimental and theoretical explanations of carbon nanotubes. Khan et al. \cite{12} inspected discharge of liquid and thermal transmission phenomena in carbon nanotube flow under the Navier-Slip condition. In single-phase and multi-phase models, Turkyilmazoglu \cite{13} explored interpretative results and thermal transport nanofluid phenomena. Manevitch et al. \cite{14} examine the impact of nonlinear radiation and optical vibrations on SWCNTS. Lee et al.
 \cite{15} looked into the heat conductivity of a fluid incorporating carbon nanotubes. Aman et al. \cite{16} investigated heat transmission enhancement in convective flow of CNTs using the Maxwell model and several types of molecular liquors. Entropy formation of Cu-water nanofluid blended convection through a cavity was inspected by Khorasanizadeh et al. \cite{17}. 
Soleimani et al. \cite{18} used core flooding to analyze the domination of a carbon nanotube-based nanofluid on oil restoration adaptability. In a critical review, Taheriang et al. \cite{19} revealed details regarding MWCNTs for improved thermophysical aspects of nanofluids. In their review study, Wang et al. \cite{20} explored the mechanics and uses of carbon nanotubes in terahertz devices. The works of \cite{21}-\cite{31} contain a number of additional major investigations with some intriguing experimental discoveries. \

As a consequence of its numerous businesslike utilizations in science and mechanics, such as agriculture, hydrology, oil and gas production, petroleum technology, civil engineering, and mining and mineral processing, researchers are currently researching the movement of fluid and transportation methodology through the exploitation of porous channels. Researchers have been spotted using Darcy's law for modelling and monitoring the discharge across porous media. For high velocities, the flow pattern becomes non-linear in a permeable medium, hence inertia effects must be considered. At incredible speeds, the Reynolds ratio (relying on porous structure) exceeds the value of unity. To compensate for such effects, Forchheimer \cite{32} united a square velocity component to the Darcian velocity in the momentum equalization. This component is known as the "Forchheimer term" by Muskat \cite{33} and it is invariably authentic for an immense Reynolds number. In Darcy-Forchheimer mixed convective movement, Seddeek \cite{34} investigated viscous dissipation and thermophoresis. A Darcy-Forchheimer outflow of magneto Maxwell liquid confined by a convectively heated sheet has been proposed by Sadiq and Hayat \cite{35}. Hayat et al. \cite{36} analyzed carbon nanotube Darcy-Forchheimer flow past a circulating disc. Iqbal et al. \cite{37} investigated carbon nanotube stagnation point flow in the companionship of a convinced magnetic field. Three-dimensional (3D) rotational leakage of water-based carbon nanotubes subjected to Darcy-Forchheimer permeable slots was studied by Hayat et al. \cite{38}. Hayat et al. \cite{39} looked into the consequences of radiant heat and dissolving thermal expansion in carbon nanotube stagnation-point discharge. The studies \cite{40}-\cite{44} provide more information on recent Darcy-Forchheimer flow research.\

The boundary layer movement of nanofluid over a shrinking/stretching plate has received much scrutiny from researchers cause of their operations in extrusions of metal sheets, extrusions of polymer sheets, and more industrial purposes. Magyari and Keller \cite{45} are the first ones to investigate the movement of heat transformation over an exponentially diminishing plate. In the literature, there are many more investigations into flow across an exponentially extending surface \cite{46}-\cite{51}.\

Suction or injection (blowing) of a fluid through the bounding surface can drastically alter the flow field. Generally, suction is known to improve skin friction, whereas injection possesses opposite effect. Fluid injection or extraction through a porous boundary wall is of fundamental relevance in boundary layer control applications such as film cooling, polymer fiber coating, and wire coating. Suction and blowing mechanisms are also essential in many technical operations, such as the design of thrust bearings and radial diffusers, as well as thermal oil recovery. Chemical processes utilize suction to remove reactants. Blowing is a technique for introducing reactants, freezing the surface, preventing corrosion or scaling, and reducing drag.\

From the above studies, the present paper addresses the flow and thermal transmission analysis of nanofluids containing carbon nanotubes (CNTS) as a nanoparticle with base fluid kerosene oil over a surface that is exponentially extending, along with an account of the phenomenon of velocity slip, thermal jump, and radiation. So yet, no such formulation has already been identified in the literature. SWCNTs and MWCNTs are the two categories of CNTs employed. The transformed non-linear governing equations are disclosed numerically by Keller-box scheme as well as bvp4c programme with the support of MATLAB software. The implication of the dimensionless criterion on velocity and thermal portraits is concluded through plots. The results of the coefficient of skin friction parameter and Nusselt number are presented in the tables.

\section{ Mathematical Interpretation of Dispute}
The standard boundary layer equations (continuity, momentum and energy) \cite{52} for this study can be interpreted as:

\begin{equation} \label{1}
\frac{\partial u}{\partial x} + \frac{\partial v}{\partial y}=0,
\end{equation}

\begin{equation} \label{2}
u\frac{\partial u}{\partial x} + v \frac{\partial u}{\partial y}=\frac{\mu_{nf}}{\rho_{nf}} \frac{\partial ^2 u}{\partial y^2} - \frac{\nu_{nf}}{K_1}u - F\ u^2-  \frac{\sigma B^2}{\rho_{nf}} u,
\end{equation}

\begin{equation} \label{3}
u \frac{\partial T}{\partial x} + v \frac{\partial T}{\partial y}=\alpha_{nf} \frac{\partial ^2 T}{\partial y^2} - \frac{1}{(\rho C_p)_{nf}} \frac{\partial q_{r}}{\partial y}, 
\end{equation}\

where $u$ and $v$ are the velocity components along $x$ and $y$ directions respectively. The $x-$axis is formed towards the order of outflow along the exponentially widening sheet, and the $y-$axis is located perpendicular to the surface. The surface is expanded and having velocity $U_w$. As demonstrated in Figure $1$, a constant magnetic field $B_0$ is injected perpendicular toward the surface. $\mu_{nf}$ refers to the dynamic viscosity of nanofluid, The density  is indicated using $\rho_{nf}$. $\nu_{nf}$ is the kinematic viscosity, $K_1$ is permeability of the porous medium and  $F$ $= \frac{c_b}{\sqrt{K_1} L}$ porous medium having a non-uniform inertia coefficient, where $c_b$ is the variety of coefficient of drag based on medium topology. $\sigma$ is  electrical conductivity, $B$ is the uniform magnetic field and $B= B_0e^\frac{x}{2L}$, where $ B_0$ is the constant magnetic field. $T$ is temperature, $\alpha_{nf}$ denotes the thermal diffusion coefficient of nanofluid where $\alpha_{nf}= \frac{k_{nf}}{(\rho C_p)_{nf}}$ with $k_{nf}$ is nanofluid temperature coefficient, ${(\rho C_p)_{nf}}$ is the specific hotness at consistent pressure and $q_{r}$ reflects the radiation heat flux of nanofluid. \

The radiative heat flux under the Rosseland assumption for radiation \cite{53}- \cite{54} can be written as:
 $$ q_{r}= -\frac{4 \sigma^*}{3K^*} \frac{\partial T^4}{\partial y} $$
 where $\sigma^*$ signifies the continual Stefan-Boltzman, $K^*$ seems to be the saturation coefficient. It is expected the discrepancy in temperature inside the discharge is relatively modest. Hence expanding $T^4$ by Taylor's series upto $T_{\infty}$ and then ignoring the terms having higher orders, then 
 $$ T^4\cong 4T^3 _\infty T - 3 T^4 _\infty, $$
 
 Now the above equation becomes:
 \begin{equation} \label{4}
 u \frac{\partial T}{\partial x} +v \frac{\partial T}{\partial y} = \alpha_{nf}\frac{\partial^2 T}{\partial y^2} + \frac{16 \sigma^* T^3_\infty }{3 (\rho C_p)_{nf} K^*} \frac{\partial^2 T}{\partial y^2}.
 \end{equation}

\begin{figure} 
\hfill
\center
{\includegraphics[width=0.9\textwidth] {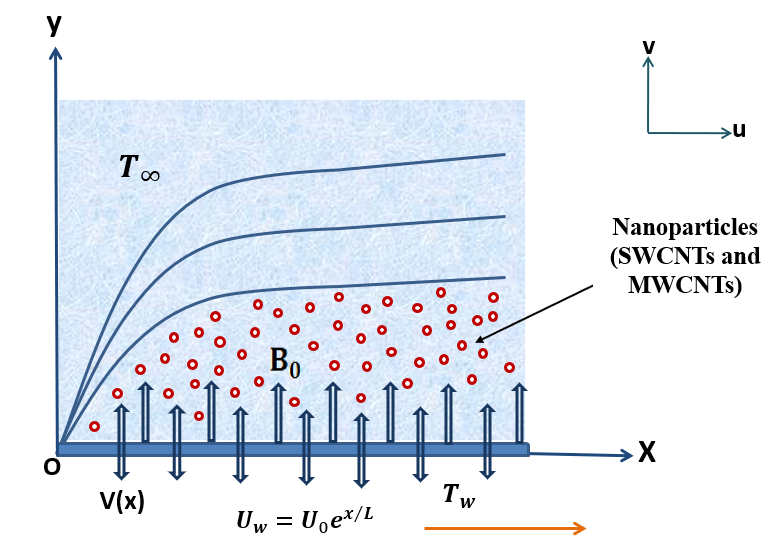}} 
\label{fig:1}
\caption{Schematic physical geometry}
\end{figure}

The thermophysical characteristics \cite{55}-\cite{58} of nanofluids are defined as:
\begin{align}     \label{5}
\begin{rcases}
\begin{aligned}
 \frac{\mu_f}{(1- \phi)^{2.5}}&= \mu_{nf}\quad\\
 \frac{\mu_{nf}}{\rho_{nf}}&= \nu{nf}, \quad\\
 (1- \phi) \rho_f + \phi \rho_{CNT}&= \rho_{nf},\quad\\
(1- \phi) (\rho C_p)_f + \phi (\rho C_p)_{CNT} &=(\rho C_p)_{nf}, \\
{K_{f}} \left( \frac{(K_{CNT}+2K_{f}) - 2 \phi (K_f -K_{CNT})}{(K_{CNT}+2K_{f}) + \phi (K_f -K_{CNT})}\right)&= {K_{nf}}.\\  
\end{aligned}
\end{rcases}
\end{align}

Here volume concentration of nanomaterial is defined by $\phi$, kinematic viscidity of nanofluid represents $\nu_{nf}$, effective viscosity of base fluid and nanofluid reflected by $\mu_f$ and $\mu_{nf}$ respectively. The basic fluid has a density $\rho_f$. The explicit warmth of base fluid and CNTs are described by $(\rho C_p)_{f}$ and $(\rho C_p)_{CNT}$ correspondingly, while potential to conduct heat of base fluid is designated through  $K_f$  and thermal conduction of CNTs are $K_{CNT}$.\

The thermophysical significance of the nanoaparticles (CNTs) and base fluid (kerosene oil) is displayed in table \ref{tab:table1}.

\begin{table}[htb]
\centering
\begin{tabular}{llll}
\hline
\textbf{Properties} & \textbf{Base Fluid}  & \multicolumn{2}{c}{\textbf{Nanoparticles}}\\
\hline
\hline
 & Kerosene Oil & SWCNT & MWCNT\\  \hline
$\rho$ & 783   & 2600  & 1600 \\
$C_p$  & 2090  & 425   & 796 \\
$K$    & 0.145 & 6600  & 3000\\ \hline
 
\end{tabular}
\caption{Data of Thermophysical properties of kerosene oil and CNTs    
\cite{59}-\cite{60}.}
\label{tab:table1}
\end{table}

\subsection{Boundary Conditions}

The corresponding boundary conditions for the study are as follows \cite{61}:
\begin{align}  \label{6}
\begin{aligned}
u&=U_w + N \nu_{nf} \frac{\partial u}{\partial y},\quad &v&=-V(x),\quad
& T&= T_w + D \frac{\partial T}{\partial y},\quad  &\mbox {at} \ y&= 0,\\
u&\rightarrow 0,\quad    &T&\rightarrow T_\infty,                \qquad    &\mbox {as} \ y &\rightarrow \infty. 
\end{aligned}
\end{align}

Here $U_w =U_0 e^{\frac{x}{L}}$ is the stretching velocity, where the reference velocity is indicated with $U_0$, $N= N_1 e^{\frac{-x}{2L}}$ is the initial velocity slip factor which changes with $x$, the minimum value of the velocity drop part is $N_1$, $T_w= T_\infty + T_0 e^\frac{x}{2L}$ is the temperature at the sheet,  $T_0$ mention  reference temperature and $D=D_1 e^{\frac{-x}{2L}}$ is the thermal slip factor which is also diverse along $x$-axis, $D_1$ is the basic value of the thermal slip factor. The no slip case will be covered if $N=D=0$. $V(x)=V_0 e^\frac{x}{2L}$ is assumed a special kind of velocity at the wall, where $V_0$ is the initial strength of suction.  $V(x)>0$ is the velocity of suction and $V(x)<0$ is the velocity of injection or blasting.

\subsection{Solution Mechanism}

The accompanying coincidence transformations are proposed to alter the aforementioned set of non-linear partial differential equations (\ref2)-(\ref4) into a set of ordinary differential equations (ODE) together with associated boundary conditions (\ref6).

\begin{align}   \label{7}
\begin{rcases}
\begin{aligned} 
u&= U_0 e^ \frac{x}{L}\ f'(\eta),\quad  &v&=-\sqrt{\frac{\nu_f U_0}{2L}} e^ \frac{x}{2L} \left[ f(\eta) + \eta f'(\eta)    \right],\\
\eta&= y\sqrt{\frac{U_0}{2\nu_f L}} e^ \frac{x}{2L},\quad     
&\psi&= \sqrt{2 \nu_f L U_0} f(\eta)\ e^ \frac{x}{2L},\\ 
T&= T_\infty + T_0 \ e^ \frac{x}{2L}  \theta(\eta).\\
\end{aligned}
\end{rcases}
\end{align}

Through injecting the stream function $\psi$, the continuity equation (\ref1) is resolved.
\begin{equation} \label{8}
\begin{split}
u=\frac{\partial \psi}{\partial y},\quad v=-\frac{\partial \psi}{\partial x}.
\end{split}
\end{equation}

The similarity variable is $\eta$, and the stream function is $\psi$.

Using the similarity transformations (\ref7), the set of translated ordinary differential equations seem to be:
\begin{align} \label{9}
\begin{rcases}  
\begin{aligned}
\frac{\mu_{nf}/ \mu_f}{\rho_{nf}/ \rho_f} f''' +ff''-2f'^2-Kf'-Fr f'^2-Mf'=0\\
\left[ \frac{1}{Pr} \frac{K_{nf}/ K_f}{(\rho Cp)_{nf}/ (\rho C_p)_f} + \frac{4}{3} R \right] \theta''-f'\theta + f \theta'=0
\end{aligned}
\end{rcases}
\end{align}
With the aid of a proper similarity transformation, boundary conditions adopt the following form:
\begin{align}   \label{10}   
\begin{aligned}
f(0)&=S,\quad  &f'(0)&=1+  A_1 \lambda f'',\quad &\theta(0)&=1+ \delta \theta',\quad     &\mbox{at}\  \eta&= 0,\\
f'(\infty)&\rightarrow 0,\quad  &\theta(\infty)&\rightarrow 0,\quad   &\mbox{as}\  \eta&\rightarrow \infty.
\end{aligned}
\end{align}

Here, $K=\frac{2 \nu_{nf} L}{K_1U_w}$ is the local porosity parameter, $F_r=\frac{c_b}{2 \sqrt{K_1}}$ is the Forchheimer number,   $M=\frac{2 \sigma B^2 L}{\rho_{f} U_w}$ is the magnetic field parameter, $R=\frac{4 \sigma^* T^3_\infty}{K_{nf} K^* }$ is thermal radiation parameter, $Pr= \frac{\nu_f (\rho Cp)_f}{K_f}$ is  Prandtl number, $S= \frac{V_0}{\sqrt{\frac{U_0 \nu_f}{2L}}}>0$ is the suction parameter ($<0$ then its blowing parameter), $\lambda= N_1 \sqrt{\frac{U_0 \nu}{2L}}$ is the velocity slip parameter and $\delta=D_1 \sqrt{\frac{U_0}{2 \nu L}}$ is the thermal slip parameter.
Also, $ A_1 =\left(  \frac{1}{(1- \phi)^{2.5} \left( 1-\phi+\phi \frac{\rho_s}{\rho_f} \right)} \right).$\

The skin friction coefficient $C_{fx}$ and local Nusselt number $Nu_x$ are defined as:
 \begin{equation}     \label{11}
\begin{split}
 C_{fx}= \frac{\tau_w}{\rho_f U_w^2},\quad     Nu_x = \frac{xq_w}{k_f(T_w-T_\infty)},\quad
\end{split}
\end{equation}\
where $\tau_w$ represents skin friction or the shear stress and $q_w$ is the surface heat flux  and is represented as:
\begin{equation}    \label{12}
\begin{split}
\tau_w = \mu_{nf} \left(  \frac{\partial u}{\partial y} \right)_{y=0} &\mbox{and}  \quad  q_w= -k_{nf} \left( \frac{\partial T}{\partial y}\right)_{y=0},
\end{split}
\end{equation}

 By using the transformation (\ref{7}) into (\ref{11}) and (\ref{12}), the skin friction number and local Nusselt number are now detailed as: 
\begin{align}     \label{13}
\begin{aligned}
\sqrt{Re}\ C_{fx} &= \frac{1}{(1- \phi)^{2.5}} f''(0),\\     
  Re_x^\frac{-1}{2}\ Nu_x &= - \frac{K_{nf}}{K_f} \sqrt{\frac{x}{2L}} \theta'(0).   
\end{aligned}
\end{align}

Here $Re_x=\frac{U_0 e^ \frac{x}{L}}{\nu_f}$ is the local Reynolds number.

\section{Computational Methodology}
\subsection{bvp4c Algorithm} 

With multiple choices of the inflow regulating parameters, the nonlinear differential equations (\ref9) and boundary conditions (\ref{10}) are quantitatively evaluated through using bvp4c simulator in MATLAB software. This three-stage Lobatto IIIa technique has been constructed using only a finite difference approach in the bvp4c solver. In order for this solver to work, users must input a set of preliminary hypotheses among the selections of boundary layer thicknesses. The findings are captured whenever the boundary conditions have been monotonically satisfied and there are not any inaccuracies in MATLAB. Before using the bvp4c algorithm, we should first represent the interconnected non-linear differential equations like a combination with first-order ODEs.

We'll simply add the predictor attributes.
\begin{align} \nonumber
\begin{aligned} 
h(1)=f, \quad h(2)=f', \quad h(3)=f''\\
h(4)=\theta, \quad  h(5)=\theta'\\
\end{aligned}
\end{align}
 
Incorporate the above different variables through into set of equations (\ref9) and generate a new network of first order equations as indicated:
\begin{align}  \nonumber
\begin{aligned}
f'&=h(2)\\
f''&=h(3)\\
f'''&= \frac{\rho_{nf}/ \rho_f}{\mu_{nf}/ \mu_f} \left[ 2h(2)^2-h(1)h(3)+Kh(2)+Fr h(2)^2+M h(2)\right] \\
\theta&=h(4)\\
\theta'&=h(5)\\
\theta''&= \left[ \left( Pr \frac{(\rho Cp)_{nf} /(\rho Cp)_f}{K_{nf}/K_f} + \frac{3}{4} R \right) \left\lbrace  h(2)h(4)-h(1)h(5) \right\rbrace  \right] \\
\end{aligned}
\end{align}

The boundary conditions (\ref{10}) are transformed into:
\begin{align}
\begin{aligned} \nonumber
h_a(1)&=S,\quad &h_a(2)&=1+A_1 \lambda h(3), \quad  h_a(4)&=1+ \delta h(5),\\ 
h_b(2)&\rightarrow 0, &\quad h_b(4)&\rightarrow 0.
\end{aligned}
\end{align}
In this case, $h_a$ identifies the location of $\eta=0$, and $h_b$ reflects the region of $\eta\rightarrow \infty$.\\

\subsection{Keller Box Scheme}
 Keller box approach is employed to calculate the governing equations (\ref9) as well as the boundary conditions (\ref{10}). The reliability of this approach is of the second series. The steps listed below can be taken to get the solution:
\begin{itemize}[noitemsep,topsep=0pt]
	\item{Reconstruct the partial differential equations toward first order ordinary differential equations.}
	\item{Employ centered-difference estimations to build finite difference interpretations of computations.}
	\item{To linearize the equations, apply Newton's linearization approach (if they are non linear).}
	\item{To tackle the restricted system, implement the block-tridiagonal-elimination method.}
\end{itemize}
\subsubsection{Transformation to First Order System}
To express equations in a first order system, we must innovate additional reliant variables $m(\eta)$, $n(\eta)$, and $o(\eta)$, as a consequence:
\begin{equation}
	f^\prime =m, \quad m^\prime=n,  \quad \theta'=o, \label{14}
\end{equation}
where $prime (\prime)$ reflects the separation with $\eta$. Accordingly, the non-linear ODEs ~(\ref{9}) can be stated as:
\begin{equation}\label{15}
\begin{aligned}
\frac{\mu_{nf}/\mu_f}{\rho_{nf}/\rho_f} n'+fn-2m^2-Km-Fr m^2-Mn&=0 ,\\ 
\end{aligned}
\end{equation}
\begin{equation}\label{16}
\begin{aligned}
\left[ \frac{1}{Pr} \frac{K_{nf}/ K_f}{(\rho Cp)_{nf}/ (\rho C_p)_f} + \frac{4}{3} R \right] o'-m \theta+fo&=0,  \\
\end{aligned}
\end{equation}

The boundary conditions ~(\ref{10}) in prospective of the newly reliant variables  are given by:
\begin{align} \label{17}
	\begin{rcases}
		\begin{aligned}
			f(0)&=S, \quad m(0)=1+\lambda n(0),  \quad \theta (0)= 1+ \delta o(0),\\
			m(\eta)&\rightarrow0, \quad \theta(\eta) \rightarrow 0. 
		\end{aligned}
	\end{rcases}
\end{align}

\subsubsection{Finite Difference Scheme}

We are now contemplating the division $\eta_{j-1} \eta_j$ along the centroid $\eta_{j-1/2}$ , which would be described as:

\begin{align}          \label{18}
	\eta_0=0, \quad \eta_j=\eta_{j-1}+h_j,   \quad \eta_j= \eta_{\infty}  
\end{align}

In which $h_j$  symbolizes the $\Delta \eta$-spacing and $\Delta \eta = 1, 2, 3,........,J$ is a progression number indicating the correspondent spot.
Now the finite difference adjustment  to equations~(\ref{14}-\ref{16}) yield the resulting sequence of equations:
\begin{align}          \label{19}
	f_{j}-f_{j-1}-h_j(m_j+m_{j-1})/2&=0,\\  
	m_{j}-m_{j-1}-h_j(n_j+n_{j-1})/2&=0,\\  
	\theta_{j}-\theta_{j-1}-h_j(o_j+o_{j-1})/2&=0,  
\end{align}

\begin{align}      
 \left( \frac{\mu_{nf}/\mu_f}{\rho_{nf}\rho_f} \right) \left( n_j-n_{j-1} \right) +\left[ \frac{h_j}{4}\right] \left(f_j+f_{j-1} \right) \left(n_j+n_{j-1} \right)-2 \left[ \frac{h_j}{4}\right] \left(m_j+m_{j-1} \right)^2 - K \left[ \frac{h_j}{2}\right]\notag\\  \left(m_j+m_{j-1} \right)  - Fr \left[ \frac{h_j}{2}\right] \left(m_j+m_{j-1} \right)^2 -M \left[ \frac{h_j}{2}\right] \left(n_j+n_{j-1} \right)&=0 \\ 
\left[ \left( \frac{1}{Pr}\right) \frac{K_{nf}/K_f}{(\rho cp) _{nf}/ (\rho cp)_f} + \frac{4}{3}R\right] \left( o_j-o_{j-1}\right) - \frac{h_j}{4} \left[ \left(m_j+m_{j-1} \right)\left(\theta_j+\theta_{j-1} \right)\right] \notag\\  + \frac{h_j}{4} \left[ \left(f_j+f_{j-1} \right)\left(o_j+o_{j-1} \right)\right]  &=0.   \label{23}            
\end{align}

Following equations are prescribed with $j=1,2,3,..... J-1$, and in addition to the boundary conditions for $j=0$ and $j=J$ are respectively:

\begin{equation} \label{24}
	\begin{rcases}
		\begin{aligned}
			f_0&=S,\quad m_0=1+ A_1\lambda n(0), \quad \theta_0= 1+ \delta o(0), \\ 
			m_J&=0, \quad \theta_J =0.
		\end{aligned}
	\end{rcases}
\end{equation}

\subsubsection{Newton's Linearization Method}

We employed Newton's method to deal with such nonlinear equations ~(\ref{19}~-~\ref{23}). For this, we've included the iterates [$f_j^i,~ m_j^i,~  n_j^i, ~ \theta_j^i, ~ o_j^i $], where $i=0,1,2,....$. We configured the appropriate for higher iterations:

\begin{align}        \label{25}
	\begin{aligned}
f_j^{i+1}&=f_j^i+\delta f_j^i,\quad &m_j^{i+1}&=m_j^i+\delta m_j^i, \quad &n_j^{i+1}&=n_j^i+\delta n_j^i, \\
\theta_j^{i+1}&=\theta_j^i+\delta \theta_j^i,\quad &o_j^{i+1}&=o_j^i+\delta o_j^i.
	\end{aligned}
	\end{align}
	
Subsequently they inserted the right-hand side of such formulations in (\ref{19}-\ref{23}), and  discarding the complex physical components within $\delta$. Generalized expressions like (remove superscript $i$ for clarification):

\begin{align}\label{26}
	\begin{rcases}
		\begin{aligned}
\delta f_{j}-\delta f_{j-1}-h_j(\delta m_j+\delta m_{j-1})/2&=(T_1)_{j},\\
\delta m_{j}-\delta m_{j-1}-h_j(\delta n_j+\delta n_{j-1})/2&=(T_2)_{j},\\
\delta \theta_{j}-\delta\theta_{j-1}-h_j(\delta o_j+\delta o_{j-1})/2&=(T_3)_{j},\\
(d_1)_j\delta n_j+(d_2)_j\delta n_{j-1}+(d_3)_j\delta f_j+(d_4)_j\delta f_{j-1}+(d_5)_j\delta m_j+(d_6)_j\delta m_{j-1}&=(T_4)_{j},\\
(e_1)_j\delta o_j+(e_2)_j\delta o_{j-1}+(e_3)_j\delta f_j+(e_4)_j\delta f_{j-1}+(e_5)_j\delta m_{j}+(e_6)_j\delta m_{j-1}\\+(e_7)_j\delta \theta_{j}+(e_8)_j\delta \theta_{j-1}&=(T_5)_{j}.
		\end{aligned}
	\end{rcases}
\end{align}
where
\begin{align*} 
(d_1)_j &= \frac{\mu_{nf}/\mu_f}{\rho_{nf}\rho_f} + \frac{h_j}{2} f_{j-1/2} , &   (d_2)_j&= (d_1)_j-2 \Bigl[ \frac{\mu_{nf}/\mu_f}{\rho_{nf}\rho_f} \Bigr]  \\
(d_3)_j&= \frac{h_j}{2} n_{j-1/2} , & (d_4)_j&=(d_3)_j,\\
(d_5)_j&= {h_j} (2+Fr) m_{j-1/2} - (K+M)  \frac{h_j}{2},& (d_6)_j&=(d_5)_j,\\
(e_1)_j & = - \left[ \left( \frac{1}{Pr}\right) \frac{K_{nf}/K_f}{(\rho cp) _{nf}/ (\rho cp)_f} + \frac{4}{3} (R) + \frac{h_j}{2} (f_{j-1/2}) \right],  \\
 (e_2)_j &=(e_1)_j - 2  \Bigl[ \left( \frac{1}{Pr}\right) \frac{K_{nf}/K_f}{(\rho cp) _{nf}/ (\rho cp)_f}  + \frac{4}{3} (R) \Bigr],\\
(e_3)_j & = \frac{h_j}{2}\left( o_{j-1/2} \right),  & (e_4)_j&=(e_3)_j,\\
(e_5)_j&= \frac{h_j}{2} h_j \left(\theta_{j-\frac{1}{2}} \right) , & (e_6)_j&=(e_5)_j,\\
(e_7)_j&= \frac{h_j}{2}\left( m_{j-1/2} \right) , & (e_8)_j&=(e_7)_j.
\end{align*}

\begin{align*}
(R_1)_j&=f_{j-1}-f_j+h_jm_{j-\frac{1}{2}},\\
(R_2)_j&=m_{j-1}-m_j+h_jn_{j-\frac{1}{2}},\\
(R_3)_j&=\theta_{j-1}-\theta_j+h_jo_{j-\frac{1}{2}},\\
(R_4)_j&=\frac{\rho_{nf}/\rho_f}{\mu_{nf}/\mu_f} (o_{j-1}-o_j)-h_j(fm)_{j-1/2} + h_j (2+Fr)m_{j-1/2}^2 + h_j (K+M)m_{j-1/2}, \\
(R_5)_j&= -\Bigl[ (Pr) \frac{(\rho cp)_{nf}/(\rho cp)_f}{(K)_{nf}/ (K)_f}  + \frac{3}{4} (R) \Bigr]( o_j-o_{j-1}) - h_j (m \theta)_{j-1/2} - h_j (fo)_{j-1/2}.
\end{align*}

For all iterates, we take the following form respectively:
\begin{align}               \label{27}
	\begin{rcases}
		\begin{aligned}
			\delta f_0 &=0, \quad \delta m_0=0, \quad \delta n_0=0,      \quad \delta \theta_0=0, \quad \delta o_0=0,\\
			\delta m_J&=0, \quad \delta \theta_J=0.
		\end{aligned}
	\end{rcases}
\end{align}
This merely indicates that the input parameters should stay unchanged across the iteration mechanism.
\subsubsection{Block Tridiagonal Structure}

The structure does have a block-tridiagonal structure, so the generalised linear pair of equations~(\ref{26}) could be resolved using the block-elimination method described by Cebeci and Bradshaw \cite{62}. Block tridiagonal structures, they are often made up of variables or constants, but in this case, an intriguing aspect is that they are made up of block matrices. The equation is transformed into a matrix vector form as follows:

\[\begin{bmatrix}
	[D_1] & [F_1]\\
	[E_2] & [D_2] & [F_2] \\
	& [E_3] & [D_3] & [F_3]\\
	& & &  &   \ddots & \cdots\\
	& & & & &   [E_{j-1}] & [D_{j-1}] & [F_{j-1}]\\
	& & &  & & & [E_j] & [D_j]
\end{bmatrix}
\begin{bmatrix}
	[\delta_1]\\
	[\delta_2]\\
	[\delta_3]\\
	\vdots\\
	[\delta_{j-1}]\\
	[\delta_J]
\end{bmatrix}
=
\begin{bmatrix}
	[T_1]\\
	[T_2]\\
	[T_3]\\
	\vdots\\
	[T_{j-1}]\\
	[T_j]
\end{bmatrix}
\]
or \begin{equation} \label{28}
	[D][\delta]=[T] 
\end{equation}
where the components of the specified heat mechanism are:

\[
\begin{bmatrix}
	\delta_1
\end{bmatrix}
=
\begin{bmatrix}
	[\delta n_0]\\
	[\delta o_0]\\
	[\delta f_1]\\
	[\delta n_1]\\
	[\delta o_1]
\end{bmatrix}
, ~~
\begin{bmatrix}
	\delta_j
\end{bmatrix}
=
\begin{bmatrix}
	[\delta m_{j-1}]\\
	[\delta \theta_{j-1}]\\
	[\delta f_j]\\
	[\delta n_j]\\
	[\delta o_j]
\end{bmatrix}
,~~ 2\le j \le J
\]

\[
\begin{bmatrix}
	T_j
\end{bmatrix}
=
\begin{bmatrix}
	(T_1)_j\\
	(T_2)_j\\
	(T_3)_j\\
	(T_4)_j\\
	(T_5)_j
\end{bmatrix}
,~~ 1\le j \le J
\]
\[
\begin{bmatrix}
	D_1	
\end{bmatrix}
= 
\begin{bmatrix}
	0 & 0 & 1 & 0 & 0 \\
	d & 0 & 0 & d & 0 \\
	0 & d & 0 & 0 & d \\
    (d_2)_1 & 0 & (d_3)_1 & (d_1)_1 & 0 \\
	0 & (e_2)_1 & (e_3)_1 & 0 & (e_1)_1 
\end{bmatrix}
\]\\

Now

\[
\begin{bmatrix}
	D_j
\end{bmatrix}
=
\begin{bmatrix}
	d & 0 & 1 & 0 & 0 \\
	-1 & 0 & 0 & d & 0 \\
	0 & -1 & 0 & 0 & d \\
    (d_6)_j & 0 & (d_3)_j & (d_1)_j & 0 \\
	(e_6)_j & (e_8)_j & (e_3)_j & 0 & (e_1)_j 
\end{bmatrix}
,~~ 2\le j \le J
\]

\[
\begin{bmatrix}
	E_j
\end{bmatrix}
=
\begin{bmatrix}
	0 & 0 & -1 & 0 & 0 \\
	0 & 0 & 0 & d & 0 \\
	0 & 0 & 0 & 0 & d \\
    0 & 0 & (d_4)_j & (d_2)_j & 0 \\
	0 & 0 & (e_4)_j & 0 & (e_2)_j 
\end{bmatrix} 
,~~ 2\le j \le J
\]

\[
\begin{bmatrix}
	F_j
\end{bmatrix}
=
\begin{bmatrix}
	d & 0 & 0 & 0 & 0 \\
	1 & 0 & 0 & 0 & 0 \\
	0 & 1 & 0 & 0 & 0 \\
    (d_5)_j & 0 & 0 & 0 & 0 \\
	(e_5)_j & (e_7)_j & 0 & 0 & 0 
\end{bmatrix} 
,~~ 1\le j < J
\]
Here $d=-\frac{h_j}{2}$.\\\\
Now we will figure out the dilemma (\ref{28}) through Lower-Upper factorisation or decomposition, supposing $D$ is any matrix whose determinant is not zero that could be decomposed into:
\begin{equation} \label{29}
	[D]=[L][U]
\end{equation}
where 
\[
\begin{bmatrix}
	L
\end{bmatrix}
=
\begin{bmatrix}
	[\alpha_1]  \\
	[\beta_2] & [\alpha_2]  \\
	& & \ddots \\
	& & \ddots & [\alpha_{J-1}]\\
	& & &  [\beta_{J}] & [\alpha_{J}] 
\end{bmatrix}
,~~
\begin{bmatrix}
	U
\end{bmatrix}
=
\begin{bmatrix}
	[I_1] & [\Gamma_1]  \\
	& [I] & [\Gamma_2]  \\
	& & \ddots &  \ddots \\
	& & & [I] & [\Gamma_{J-1}]\\
	& & & & [I] 
\end{bmatrix}
\]
While [$I$] represents the $5 \times 5$ matrix as well as [$\alpha j$], [$\Gamma j$], and [$\beta j$] are $5 \times 5$ ordered matrices for whom components are predicted through consecutive expression:

\begin{align}       \label{30}
\begin{rcases}
[\alpha_1] & =[D_1]\\
[D_1][\Gamma_1]&=[F_1]\\
[\alpha_j]&=[D_j]-[\beta_j][\Gamma_{j-1}],~~ j=2,3,...J \\
[D_j][\Gamma_j]&=[F_j],~~ j=2,3,...J
\end{rcases}
\end{align}
By replacing  equation~(\ref{29}) in~(\ref{28}), we acquire:
\begin{equation}   \label{31}
	[L][U][\delta]=[R]
\end{equation}
Now allow us to define,
 \begin{equation} \label{32}
	[U][\delta]=[G]
\end{equation}
As just a conclusion, equation~(\ref{31}) produces
\begin{equation} \label{33}
	[L][G]=[R], 
\end{equation}
where \[
\begin{bmatrix}
	G
\end{bmatrix}
=
\begin{bmatrix}
	[G_1]\\
	[G_2]\\
	\vdots\\
	[G_{j-1}]\\
	[G_j]
\end{bmatrix}
,~~ 1 \le j \le J
\]

In this case, [$G j$] are column matrices of dimension $5 \times 1$, and entries can be acquired by mathematical analysis (\ref{33}) like:

\begin{equation}
	[\alpha_1][G_1]=[T_1], ~~~
	[\alpha_j][G_j]=[T_j]-[\beta_j][G_{j-1}],~~ 2 \le j \le J.
\end{equation}
The elements of $\Gamma j$, $\alpha j$, and $G j$ are determined adopting forward sweep, while the phonemes of $\delta$ are accurately assessed from the equation (\ref{32})  utilizing backward sweep, and thus the constituents include:
\begin{equation}
	[\delta_J]=[G_J], ~~~
	[\delta_j]=[G_j]-[\Gamma_j][\delta_{j+1}],~~ 1 \le j \le J-1.
\end{equation}
Such iterative techniques are carried out till the acceptable convergence criterion is accomplished i.e. $  |{\delta v_0}^{(i)}| < \epsilon $, where $\epsilon $ is the intended level of precision.

\section{Results and Discussion}
\subsection{Code Verification and Parameter Selection}
In the present study, the transmission of heat and circulation of CNTs (SWCNT and MWCNT) with lamp fuel oil as just a base fluid under the influence of thermal jump and velocity slip across an exponentially elaborated plate were explored through the proposed investigation. First, the asymmetric partial differential equations (\ref 2)-(\ref 4) are modified into a group of nonlinear ordinary differential equations (\ref 9) through exploiting appropriate similarity transformations. Then, the set of non-linear differential equations is numerically computed by Keller-box scheme and bvp4c algorithm in MATLAB. Table $2$ demonstrates the comparative analysis for the character of rate of thermal conduction $\theta'(0)$ in the non-appearance of magnetic field i.e. $M=0$, nanoparticle concentrations i.e. $\phi=0 $, Forchheimer parameter $Fr=0$, porosity parameter $K=0$ and the radiation effect i.e. $R=0$ for $Cu-$ water based nanofluid for the divergent measurements of $Pr$. The acquired consequences are correlated with Magyari and Keller \cite{45}, Ishak \cite{63} and EI. Aziz \cite{64}. They have been established in compliance with historically published studies, as demonstrated in table $2$. Plots depict the effects of important aspects on profiles of velocity and temperature. Also, the effect of skin friction coefficient parameter and Nusselt number are portrayed through table $3$ and $4$ correspondingly. By retaining other factors constant, i.e. $\phi=0.2, Fr=0.25, K=1, M=2, S=0.1, \lambda=0.1, \delta=0.1$, the effect of various flow controlling parameters on the velocity and thermal sketches is demonstrated in figures $2$ and $3$ correspondingly.
\begin{table}[h!]
\centering
\begin{tabular}{|l|l|l|l|l|l|l|}
\hline
\textbf{Pr} & \textbf{Magyari and Keller \cite{45}}  & \textbf{Ishak \cite{63} } & \textbf{EI Aziz \cite{64} }  & \textbf{Present study}\\
\hline
 1  & 0.954782    & 0.9548 & 0.954785   & 0.955885\\
 2  &  ----       &  ----  & ----       & 1.47122\\
 3  & 1.869075    & 1.8691 & 1.869074   & 1.868878\\
 5  & 2.500135    & 2.5001 & 2.500132   & 2.499982\\
 10 & 3.660379    & 3.6604 & 3.660372   & 3.660255\\ \hline
\end{tabular}
\caption{Comparison to  the past studies for temperature values, i.e. $-\theta'(0)$ throughout the exclusion of nanofluid, where $K=Fr=0$.}
\label{Table:ta}
\end{table}

\subsection{Graphical and Tabular Analysis}
\begin{figure}[htb] 
\subfigure[Variation in $f'(\eta)$ with $\phi$ against $\eta$]
{\includegraphics[width=0.46\textwidth] {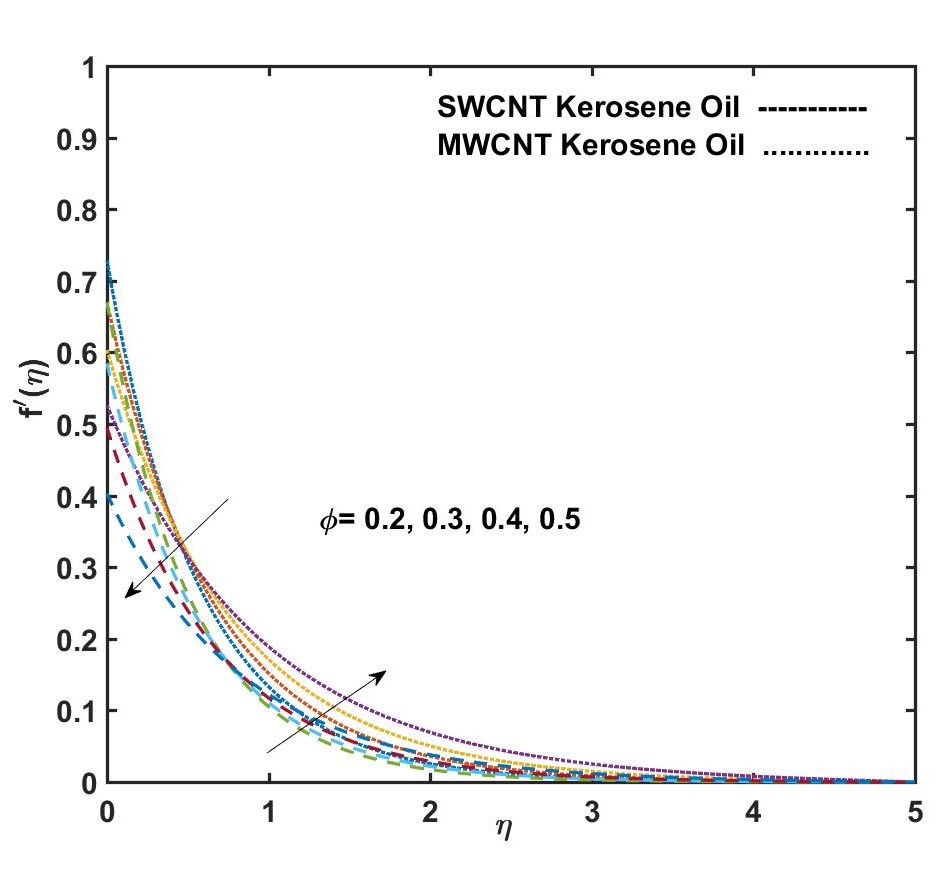}} 
\label{fig:1}
\hfill
\subfigure[Variation in $f'(\eta)$ with $Fr$ against $\eta$]
{\includegraphics[width=0.53\textwidth] {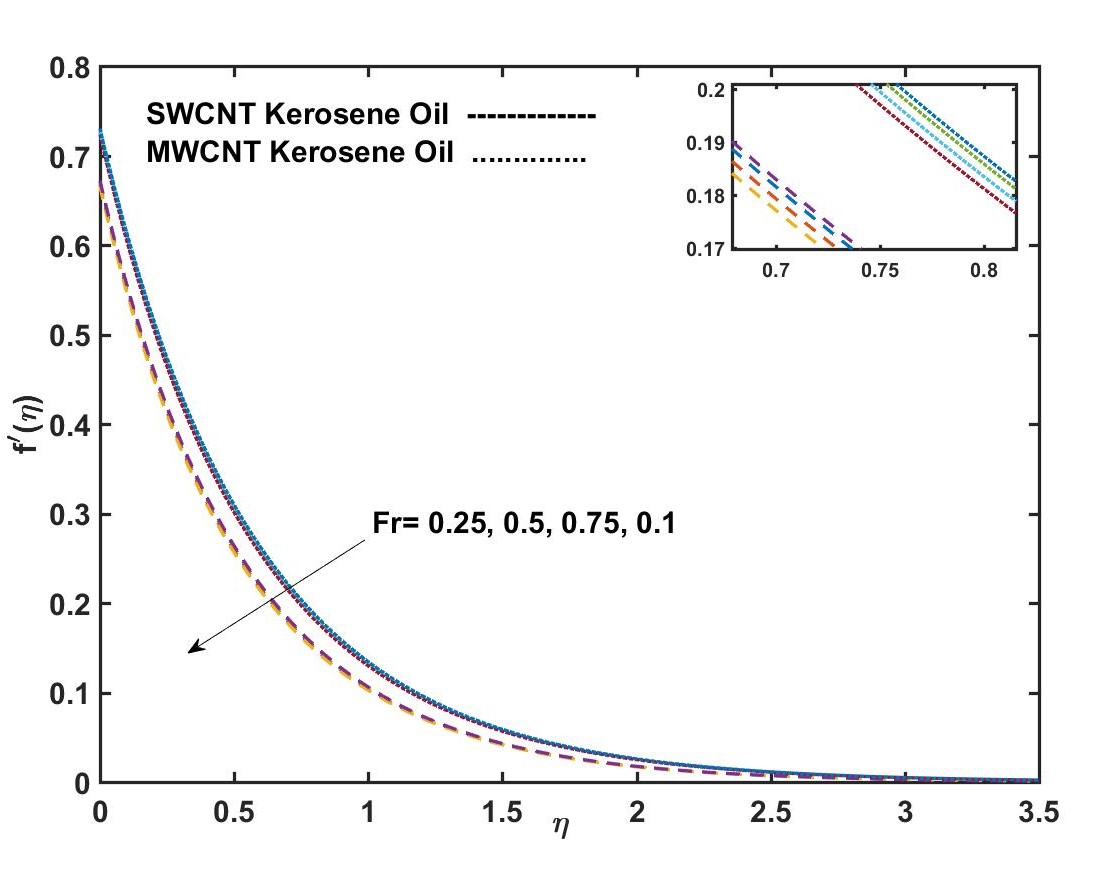}} 
\label{fig:2}
\hfill
\end{figure}

\begin{figure} [htb]
\subfigure[Variation in $f'(\eta)$ with $K$ against $\eta$]
{\includegraphics[width=0.50\textwidth] {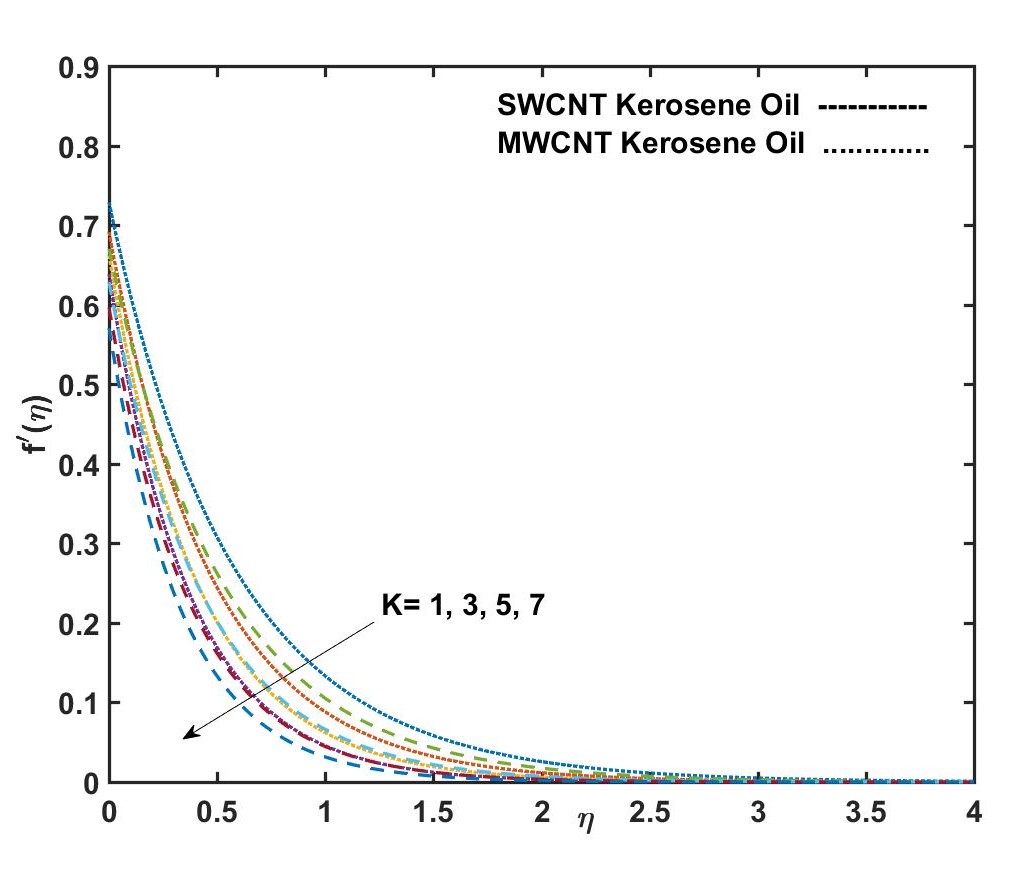}} 
\label{fig:3}
\hfill
\subfigure[Variation in $f'(\eta)$ with $M$ against $\eta$ ]{\includegraphics[width=0.48\textwidth] {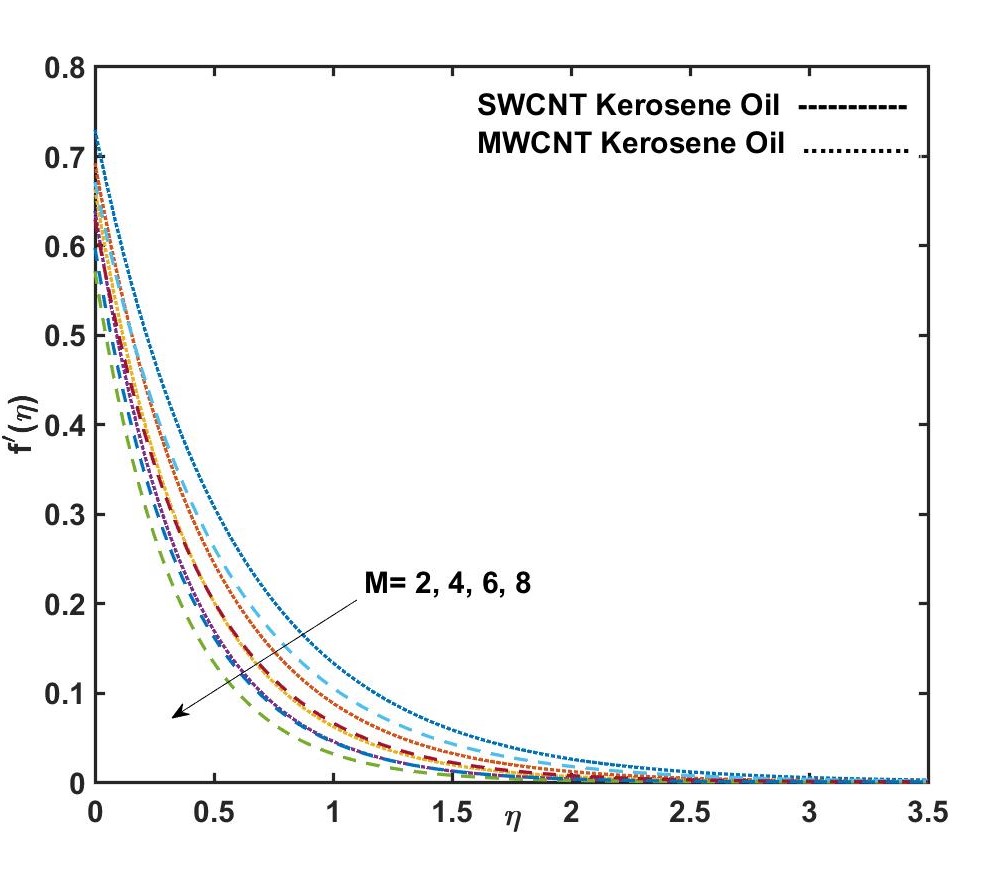}} 
\label{fig:4}
\hfill
\subfigure[Variation in $f'(\eta)$ with $S$ against $\eta$ ]{\includegraphics[width=0.47\textwidth] {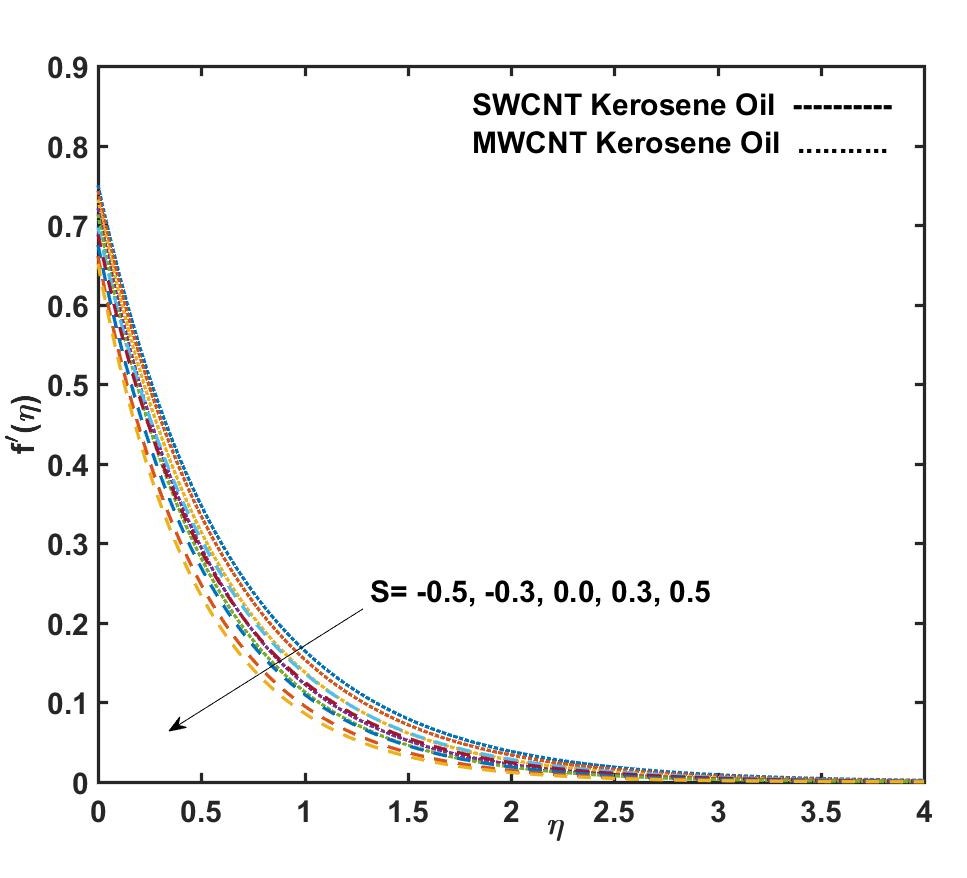}} 
\hfill
\subfigure[Variation in $f'(\eta)$ with $\lambda$ against $\eta$]{\includegraphics[width=0.47\textwidth] {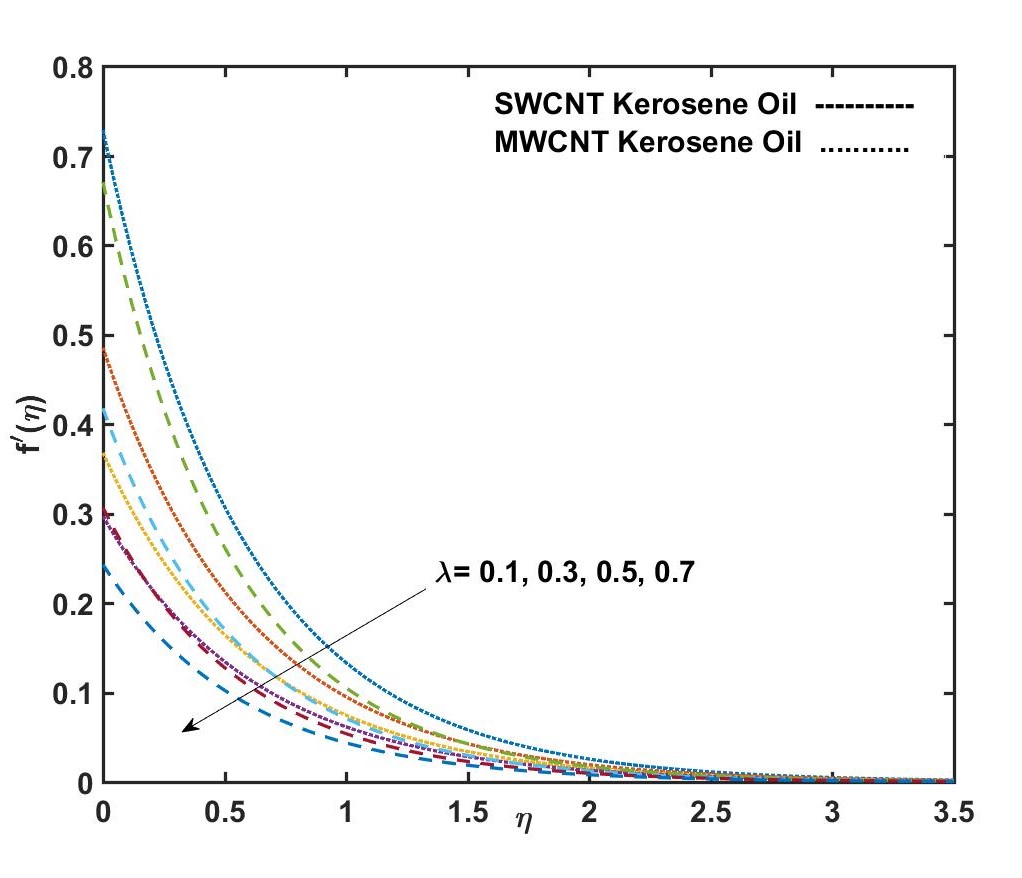}} 
\label{fig:4}
\hfill
\label{fig:5}
\caption{Alteration in $f'(\eta)$ for the distinct values of $\phi$, $Fr$, $K$, $M$, $S$ and $\lambda$}
\end{figure}
\subsubsection{Variation in Velocity Profiles}
The velocity profiles of  $f'(\eta)$ over an exponentially stretched sheet is examined in figure $2$ by varying flow regulating parameters. Inside the confined zone, the figure $2(a)$ demonstrates dual phenomena with mounting $\phi$ values.  For SWCNTs, the region $0\leq \eta < 0.6$, shows decreasing behaviour, while at $0.6\leq \eta \leq 1.0$, it shows increasing demeanor with the soaring values of $\phi$.  For MWCNTs, in the vicinity $0\leq \eta < 0.5$, it commit collapsing nature, whereas the area $0.5\leq \eta \leq 1.0$, with increasing values of $\phi $, it displays increasing category.  The velocity distribution is found to be the leading function for the enhancement in $\phi$ for SWCNTs and MWCNTs. The accession in the percentage of nanomaterials causes an increase in the convective flow. Also, it has  been observed that $f'(\eta)$ increases for lamp fuel oil nanoliquid as long as MWCNTs when related to SWCNTs.\

Figure $2(b)$ is outlined to inspect that, by virtue of the velocity profile $f '(\eta)$ is distressed by the accumulation in $Fr$. From the sketch, it has been illustrated that the velocity graph reduces for both SWCNTs and MWCNTs with the increasing values of the inertia coefficient $Fr$. 

\begin{figure}[h!] 

\subfigure[Variation in $\theta(\eta)$ with $\phi$ against $\eta$]{\includegraphics[width=0.44\textwidth] {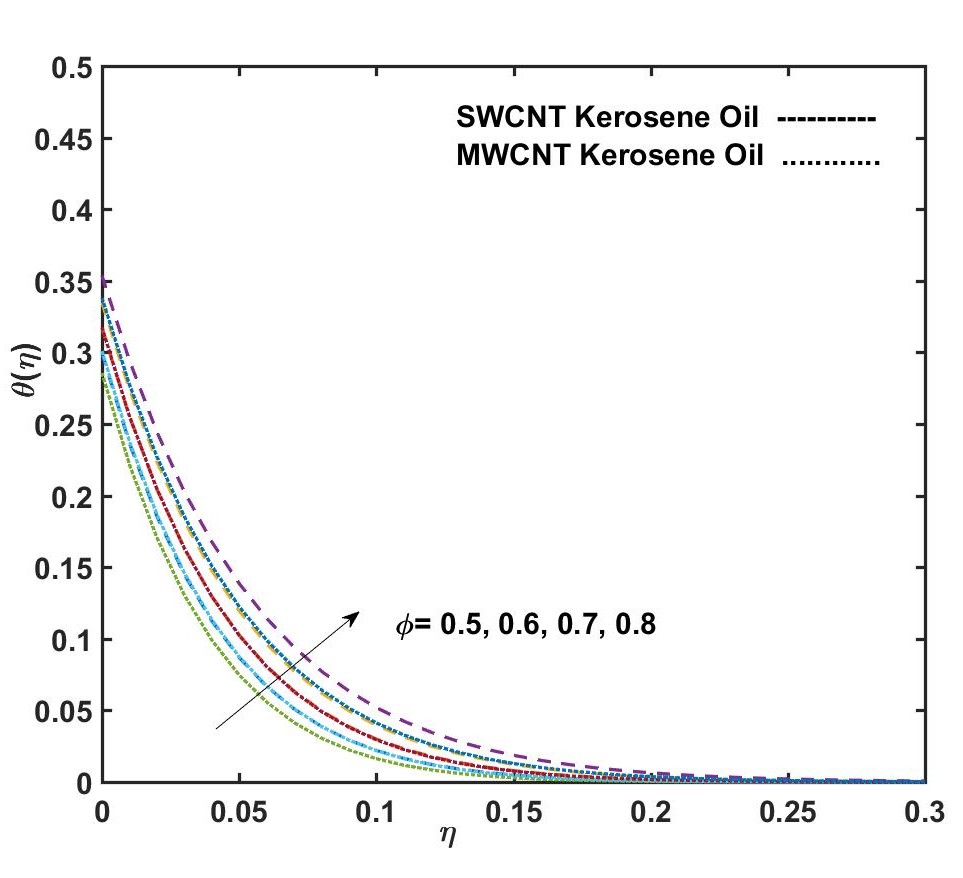}} 
\label{fig:6}
\hfill
\subfigure[Variation in $\theta(\eta)$ with $K$ against $\eta$]{\includegraphics[width=0.45\textwidth] {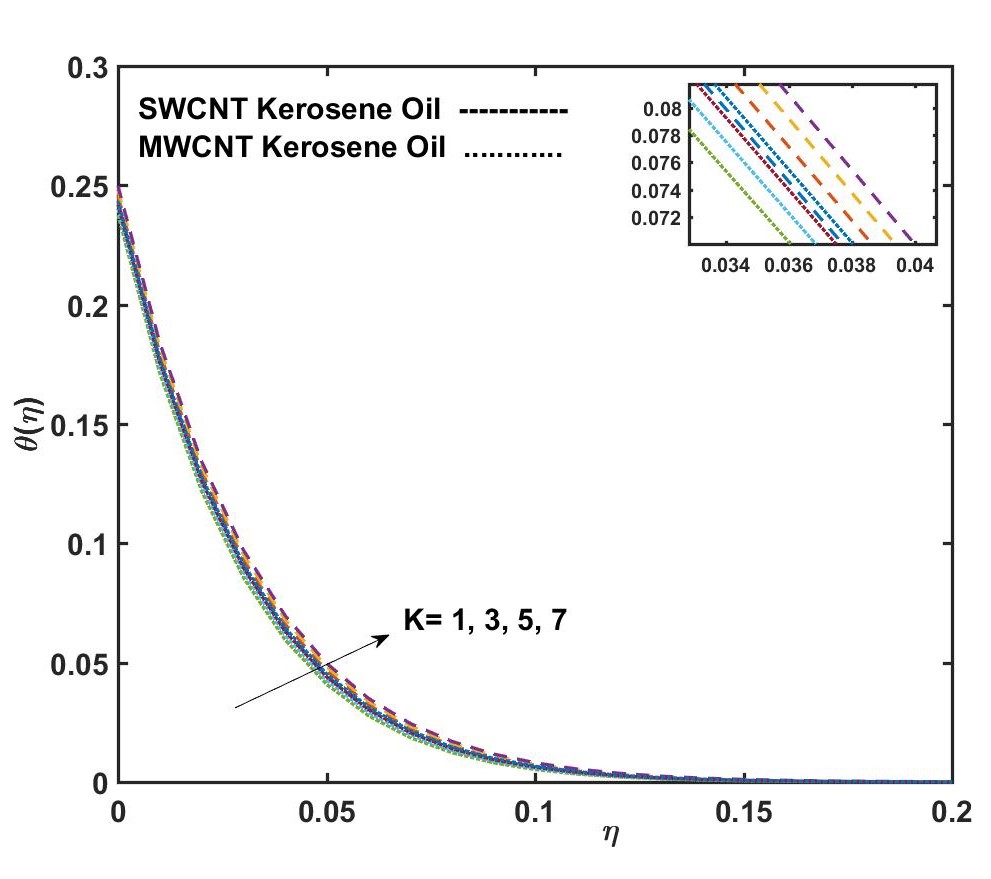}} 
\label{fig:7}
\end{figure}

\begin{figure}[h!] 
\subfigure[Variation in $\theta(\eta)$ with $M$ against $\eta$]{\includegraphics[width=0.47\textwidth] {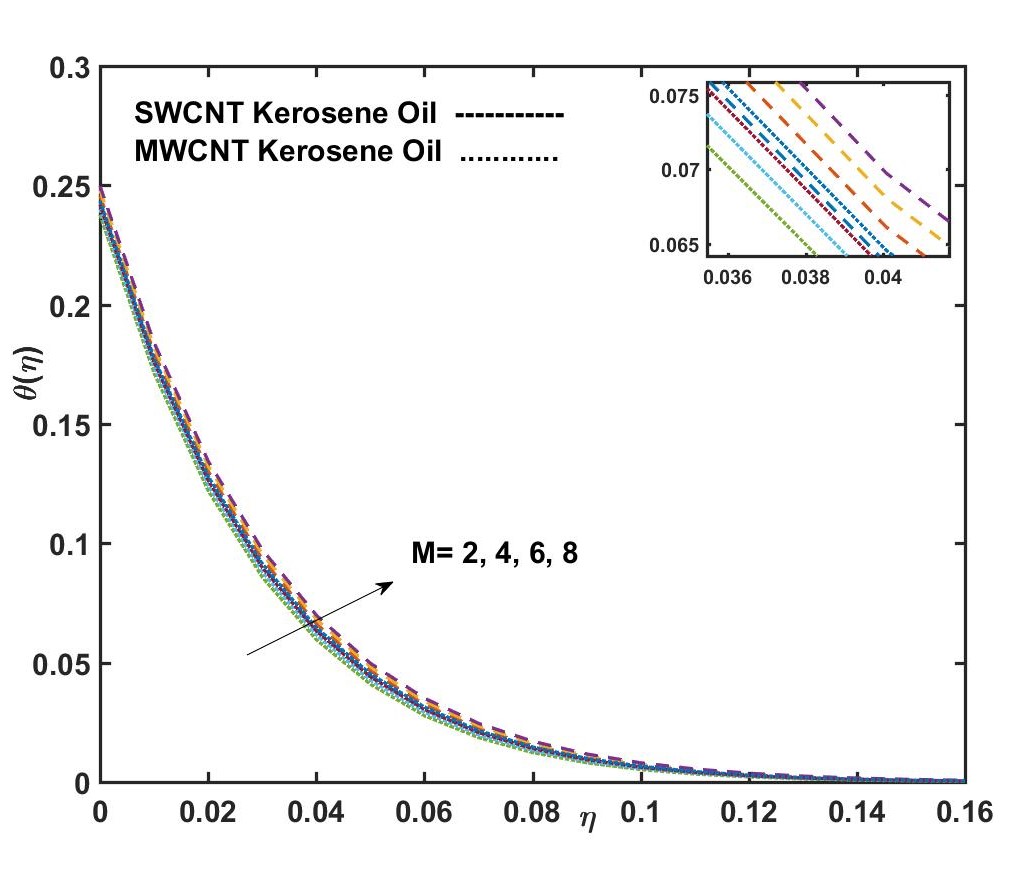}} 
\label{fig:8}
\hfill
\subfigure[Variation in $\theta(\eta)$ with $S$ against $\eta$]{\includegraphics[width=0.49\textwidth] {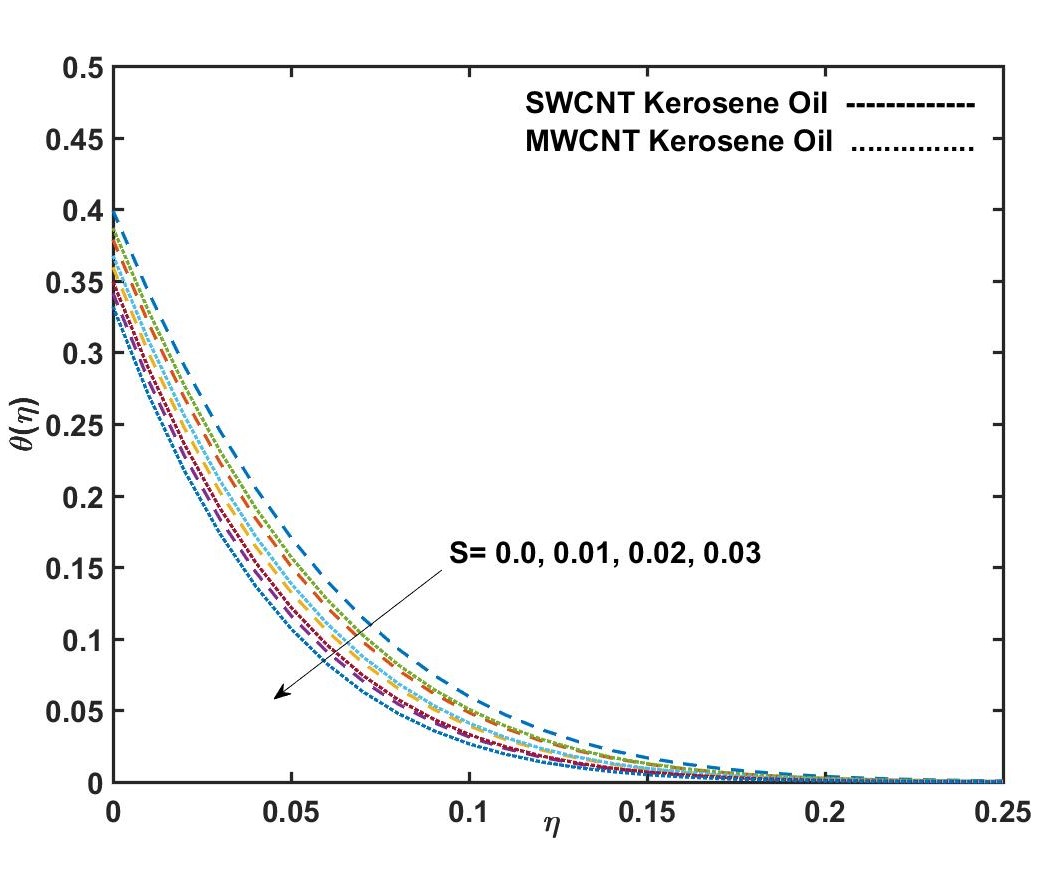}} 
\label{fig:9}
\hfill
\subfigure[Variation in $\theta(\eta)$ with $\lambda$ against $\eta$] {\includegraphics[width=0.47\textwidth] {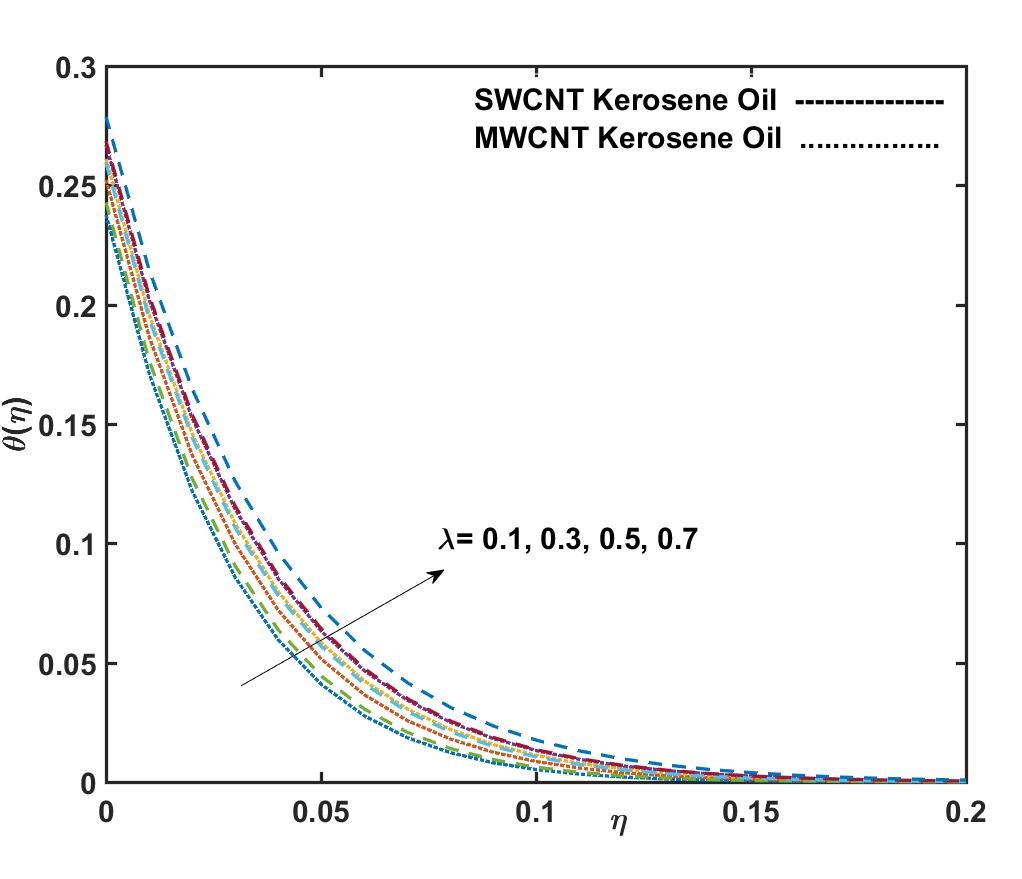}} 
\label{fig:10}
\hfill
\subfigure[Variation in $\theta(\eta)$ with $\delta$ against $\eta$] {\includegraphics[width=0.49\textwidth] {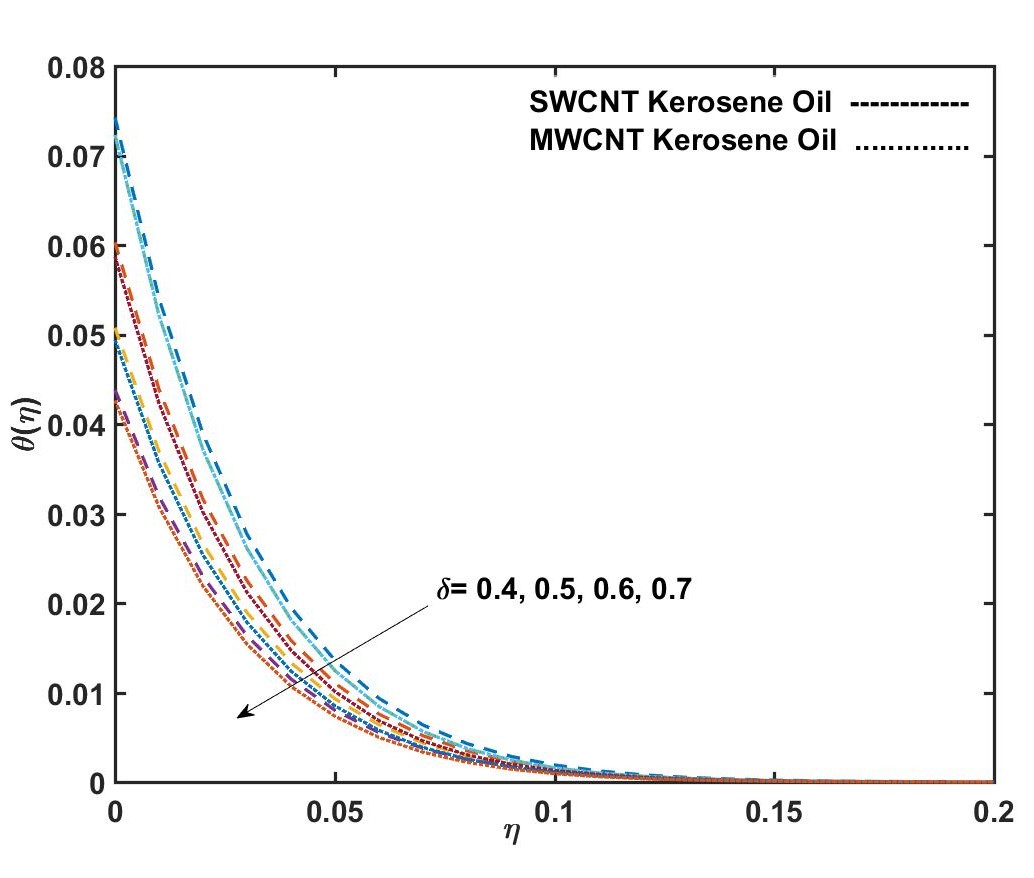}} 
\label{fig:11}
\end{figure}
\pagebreak
\begin{figure} 
\hfill
\center
\subfigure[Variation in $\theta(\eta)$ with $R$ against $\eta$]{\includegraphics[width=0.44\textwidth] {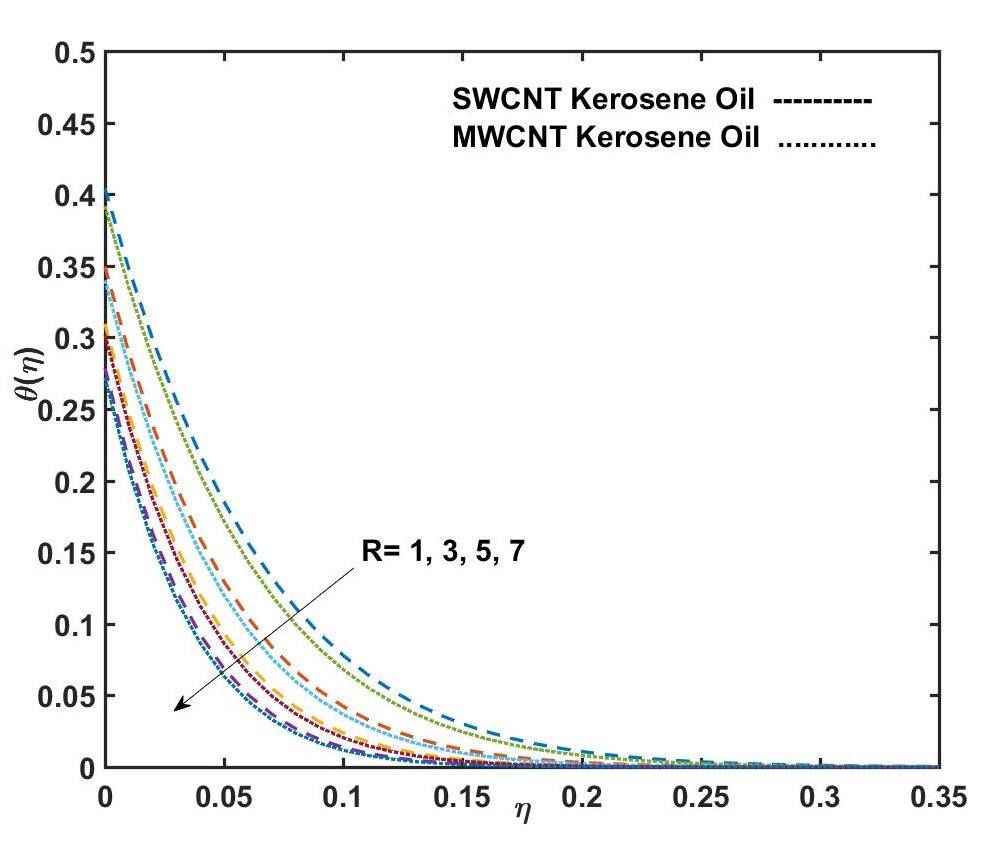}} 
\label{fig:12}
\caption{Alteration in $\theta(\eta)$ for the distinct values of $\phi$, $K$, $M$, $S$, $\lambda$, $\delta$ and $R$}
\end{figure}

 Figure $2(c)$ interprets the change in velocity field $f'(\eta)$ with the fluctuating values of the porosity parameter $K$. It has been prove that with the larger values of $K$, the velocity profile $f'(\eta)$ decreases simultaneously for both SWCNTs and MWCNTs. For lamp fuel oil based nanofluids, the flow appears to significantly slow down in the axial direction at high values of the porosity parameter $K$. Apparently, the flow intransigence hikes whenever a porous medium is introduced. Furthermore, for lamp fuel oil base liquid, the velocity distribution of SWNTs is wider than that of MWNTs. \

 Figure $2(d)$ highlights the implication of the boosting values of magnetic parameter $M$ on the velocity field. It is concluded that for the increment in the attitude of the magnetic field, the velocity profile shows a reduction in their behaviour for SWCNTs and MWCNTs. For higher magnetic estimator predictions, the velocity $f'(\eta)$ and associated boundary layer diminish simultaneously. The upsurge in $M$ represents a spike in insusceptible potential (Lorentz force), and hence the velocity of the liquid decelerates. It is also reported that MWCNT kerosene oil seems to have a dominating velocity distribution whenever compared to SWCNT kerosene oil.\
 
  Figure $2(e)$ depicts the dominance of the suction or blowing parameter $S$ on the velocity profile $f'(\eta)$. From the figure, it has been detected that the velocity field downturns significantly with the increment in the suction parameter $S$ whereas the blowing parameter shows an increment in the velocity field. Whenever the wall suction $(S > 0)$ is addressed, the boundary layer density diminishes and the velocity field reduces. $S=0$ reflects the scenario of a non-porous stretched plate. It is discovered that raising the suction parameter greatly decreases fluid velocity, whereas blowing increases the fluid velocity. \

The effect of velocity slip parameter $\lambda$ on velocity profile $f'(\eta)$ is delineated in figure $2(f)$. The velocity graph $2(f)$ illustrates that as the sheet's distance $(\eta)$ climbs, the transit rate declines. The velocity disappears at a certain large distance from the sheet (at $\eta=3.5$ in all cases). As a result, horizontal velocity decreases as $\lambda$ increases. When the layer slides, the flow rate all around the plate never approaches the elongating mobility of the material. Because the stretching sheet's pulling can only be moderately disclosed under the slide condition to the liquid, similar slip velocity increases as $\lambda$ increases. It should be noted that $\lambda$ has a momentous impact on the CNTs. Also, it is recognised that the kersoene-SWCNTs have higher velocity than kerosene MWCNTs.\ 

\subsubsection{Variation in Thermal Profiles}

Figure $3$ scrutinize the consequences of diverse outflow predominating parameters on the thermal profile  $\theta (\eta)$ over an exponentially stretched sheet. The ramification of expanding values of volume fraction $\phi$ on the thermal profile $\theta (\eta)$ is portrayed in figure $3(a)$. It is noticed that the temperature profile increases with the increasing character of $\phi$ for both SWCNTs and MWCNTs. This is caused by the addition of carbon nanotubes. The thermal conductivity of the fluid increases vastly because of the very high rate of thermal transmission of CNTs in comparison with other nanomaterials. Moreover, it is also noticed that the temperature profile for SWCNTs is greater as compared to MWCNTs. \

Figure $3(b)$ depicts the change in temperature portrait as the porosity parameter $K$ is varied. For the larger values of $K$, an enhancement is noticed in the temperature profile $\theta(\eta)$ for both SWCNTs and MWCNTs. The temperature profile in figure $3(c)$ is seen to be significantly increased with the escalating values of magnetic parameter $M$ for both SWCNTs and MWCNTs. The existence of a magnetic field $M$ reduces the density of the momentum boundary layer while increasing the density of the thermal boundary layer. \

Figure $3(d)$ depicts how a change in the suction/blowing parameter $S$ alters the thermal sketch $\theta(\eta)$, resulting a reduction in the nature of both SWCNTs and MWCNTs. From the figure, it is seen that the temperature profile decelerates as the suction variable $S$ is elevated. As a result, the effect of the suction/blowing parameter is to maximise the amount of thermal transmission from the surface to the neighbouring fluid. Bringing vacuum fluids virtually towards the interface region means that solutions exhibit a spacious flow separation, and thus the density of a momentum pressure gradient is reduced.\

The behaviour of $\theta(\eta)$ for varying ethics of the velocity lapse parameter $\lambda$ is advertised in figure $3(e)$. It is found that the thermal sketch shows that with the increment in the character of the velocity blunder framework $\lambda$, the temperature profile also boosted for both SWCNTs and MWCNTs. From the figure, it has been detected that the thermal portrait is greater for SWCNTs as correlated to MWCNTs. \

Meanwhile, temperature jump $\delta$ in figure $3(f)$, shows a reduction with the increment in $\delta$ for both SWCNTs and MWCNTs. The temperature profile initially declines with $\delta$, but over a certain separation $\eta$ from the sheet, this component fades entirely. As the  $\delta$ improves, limited warmth is transported towards the liquid from the plate, and therefore the temperature decreases gradually. Figure $3(g)$ depicts the behaviour of the temperature profile $\theta(\eta)$ with the ascending radiation parameter $R$. It has been established that when $R$ is raised, the surface heat transmission capacity comes to a halt.\

\begin{table}[h!]
\resizebox{1.11\textwidth}{!}{\begin{minipage}{\textwidth}
\center
\begin{tabular} {lllllllllccc}
\hline
 {$\phi$} &  {$M$} & {$Fr$}  & {$K$}  & $S$ &  $\lambda$ &   \multicolumn{2}{c}{ {$-\frac{1}{(1-\phi)^{2.5}} f''(0)$}} \\
   &    &   &  &   &      & SWCNT  & MWCNT & \\
\hline
0.02  & 0.1 & 0.2 & 0.1   & 0.1 & 0.1      & 1.2486  & 1.2408  & \\ 
0.05  &     &     &          &     &       & 1.3070  & 1.2894  &\\
0.07  &     &     &          &     &       & 1.3455  & 1.3227  &\\
0.1   & 1.0 &     &          &     &       & 1.6504  & 1.6211  &\\
      & 1.5 &     &          &     &       & 1.7640  & 1.7352  &\\
      & 2.0 &     &          &     &       & 1.8658  & 1.8379  &\\ 
      & 2.5 & 0.1 &          &     &       & 1.9493  & 1.9222  &\\        
      &     & 0.2 &          &     &       & 1.9581  & 1.9315  &\\
      &     & 0.3 &          &     &       & 1.9668  & 1.9406  &\\  
      &     & 0.4 & 0.1      &     &       & 1.9754  & 1.9497  &\\   
      &     &     & 0.3      &     &       & 2.0095  & 1.9843  &\\             
      &     &     & 0.5      &     &       & 2.0424  & 2.0179  &\\                 
      &     &     & 0.7      & 0.1 &       & 2.0744  & 2.0504  &\\                                                      
      &     &     &          & 0.2 &       & 2.1074  & 2.0824  &\\                                     
      &     &     &          & 0.4 &       & 2.1746  & 2.1476  &\\                                         
      &     &     &          & 0.5 & 0.1   & 2.2088  & 2.1807  &\\                                             
      &     &     &          &     & 0.4   & 1.1861  & 1.2330  &\\                                                 
      &     &     &          &     & 0.7   & 0.8184  & 0.8678  &\\                                                    
      &     &     &          &     & 0.9   & 0.6793  & 0.7260  &\\                                                        

\hline
\end{tabular}
\end{minipage}}\\
\caption{Alteration in {$-\frac{1}{(1-\phi)^{2.5}} f''(0)$} for a variety of input variables with $Pr=21$.}
\label{Table:ta}
\end{table}

\begin{table}[h!]
\resizebox{1.11\textwidth}{!}{\begin{minipage}{\textwidth}
\center
\begin{tabular} {lllllllllccc}
\hline
 {$\phi$} &  {$M$} & {$Fr$}  & {$K$} & {$R$}  & $S$ &  $\lambda$ &  $\delta$ &   \multicolumn{2}{c}{{$-{\frac{K_{nf}}{K_f}} \theta'(0)$}}\\
   &    &   &  &   &   &  &   & SWCNT  & MWCNT & \\
\hline
0.02  & 0.1 & 0.2 & 0.1 & 10  & 0.1 & 0.1 & 0.1  & 8.5023   & 8.5062   & \\ 
0.05  &     &     &     &     &     &     &      & 9.1959   & 9.2067   &\\
0.07  &     &     &     &     &     &     &      & 9.6799   & 9.6963   &\\
0.1   & 1.0 &     &     &     &     &     &      & 10.4209  & 10.4493  &\\
      & 1.5 &     &     &     &     &     &      & 10.4110  & 10.4406  &\\
      & 2.0 &     &     &     &     &     &      & 10.4020  & 10.4327  &\\ 
      & 2.5 & 0.1 &     &     &     &     &      & 10.3946  & 10.4261  &\\        
      &     & 0.2 &     &     &     &     &      & 10.3938  & 10.4253  &\\
      &     & 0.3 &     &     &     &     &      & 10.3930  & 10.4246  &\\  
      &     & 0.4 & 0.1 &     &     &     &      & 10.3923  & 10.4239  &\\   
      &     &     & 0.3 &     &     &     &      & 10.3891  & 10.4212  &\\             
      &     &     & 0.5 &     &     &     &      & 10.3861  & 10.4184  &\\                 
      &     &     & 0.7 & 1   &     &     &      & 8.7329   & 8.8044   &\\                     
      &     &     &     & 5   &     &     &      & 9.6695   & 9.7172   &\\                         
      &     &     &     & 7   &     &     &      & 9.9964   & 10.0369  &\\
      &     &     &     & 10  & 0.1 &     &      & 10.3832  & 10.4158  &\\                                 
      &     &     &     &     & 0.2 &     &      & 11.3612  & 11.3788  &\\                                     
      &     &     &     &     & 0.4 &     &      & 12.1877  & 12.1964  &\\                                         
      &     &     &     &     & 0.5 & 0.1 &      & 12.3909  & 12.3978  &\\                                             
      &     &     &     &     &     & 0.4 &      & 12.3830  & 12.3903  &\\                                                 
      &     &     &     &     &     & 0.7 &      & 12.3799  & 12.3871  &\\                                                    
      &     &     &     &     &     & 0.9 & 0.1  & 12.3787  & 12.3858  &\\                                                        
      &     &     &     &     &     &     & 0.15 & 8.4542   & 8.4574   &\\                                                        
      &     &     &     &     &     &     & 0.17 & 7.5027   & 7.5053   &\\                                                            
      &     &     &     &     &     &     & 0.19 & 6.7438   & 6.7458   &\\                                                        
\hline
\end{tabular}
\end{minipage}}\\
\caption{Alteration in {$-{\frac{K_{nf}}{K_f}} \theta'(0)$} for a variety of input variables with $Pr=21$.}
\label{Table:ta}
\end{table}

\subsubsection{Variation in Coefficient of Skin Friction and Rate of Heat Transmission}

The divergence in skin friction coefficient and Nusselt number estimates with the growing ethics of volume fraction $\phi$, magnetic parameter $M$, Forchheimer parameter $Fr$, porosity parameter $K$, radiation parameter $R$, suction parameter $S$, velocity slip parameter $\lambda $ and temperature jump parameter $\delta$ are displayed in table $3$ and $4$. In table $3$, an enhancement in the SWCNTs and MWCNTs is made up of the surged volume fraction of the nano materials, magnetic number, Forchheimer number, porosity parameter, and suction parameter. However, SWCNT and MWCNT were reduced with the improving values of $\lambda$. However, the coefficient of skin friction and rate of heat transmission climb at the surface as the nanomaterial volume fraction rises. This is because adding CNTs to nanofluids makes them more viscous, which improves heat conductivity. From the table $4$, it is seen that the Nusselt number inflates in the pair of SWCNTs and MWCNTs with the accession values of $R$ and $S$ whereas it decreases with the increment in $M$, $Fr$, $K$, $\lambda$ and $\delta$.

\section{Conclusion}

In the current investigation, the 2-dimensional boundary layer ciculation of nanofluid containing carbon nanotubes (CNTs) under the influence of thermal radiation, velocity slip, and temperature fluctuation throughout an exponentially extending layer is addressed. In this work, lamp fuel oil is employed as a conventional fluid for the extraction of CNT-based nanoparticles. There are two standard forms of CNTs that were introduced, named single-walled carbon nanotubes (SWCNTs) and multi-walled carbon nanotubes (MWCNTs) for concentration in the working fluid. Both the bvp4c algorithm and the Keller-box technique in MATLAB software are employed to evaluate the set of nonlinear equations (\ref {9}) and boundary conditions (\ref{10}). The consequences of various flow modulating non-dimensional characteristics on velocity, temperature, coefficient of skin friction , and rate of thermal transmission are tested and provided graphically and tabularly. The key themes of the ongoing investigation are listed below:
\begin{itemize}

\item{When the nanoparticle concentration climbs, MWCNTs generate an inflate in the velocity field relative to SWCNTs and a dual solution occurs in the velocity field. }

\item{The velocity diminishes as the velocity slip grows. Temperature is reduced as the thermal blunder framework is improved.}

\item{The suction criterion $S$ dwindled the fluid velocity and temperature with their growing values. Also, the temperature field is reduced as the radiation parameter $R$ is enhanced. }
 
\item{An enhancement is engraved in skin friction coefficient with the increasing values of volume fraction of nanoparticles $\phi$, magnetic parameter $M$, inertia coefficient $Fr$, porosity parameter $K$ and suction $S$ for both considered nanoparticles (SWCNTs and MWCNTs), while it reduces for the pair of (SWCNTs and MWCNTs) with the increment in velocity slip. }

\item{ The heat transmission rate shows an increment with the improvement in volume fraction of nanoparticles, radiation and suction parameters respectively. However, an opposite trend is examined in magnetic parameters, inertia coefficient, porosity parameter, velocity slip, and temperature jump for SWCNTs and MWCNTs.}

\end{itemize}

\section*{Compliance with Ethical Standards}

\subsection*{Declarations of Originality}
The author declares that no similar work has been done earlier by anyone and takes full responsibility of the originality of the figures and graphs presented in this work.

\subsection*{Funding}
Not applicable.

\subsection*{Conflicts of interest}
There is no conflict of interest.

\subsection*{Availability of data and material}
Not applicable.

\end{titlepage}
\end{document}